\documentclass[a4paper, 12pt, draft]{amsart}
\usepackage{graphicx}
\usepackage{verbatim}
\usepackage[noadjust]{cite}
\newtheorem*{theorem}{Theorem}%[section]

\theoremstyle{definition}

\newtheorem*{definition}{Definition}

\theoremstyle{remark}

\numberwithin{equation}{section}
\def\N{{\mathbb N}}
\def\R{{\mathbb R}}

\newcommand{\Xs}{X^{\ast}}
\newcommand{\xs}{x^{\ast}}

\newcommand{\ys}{y^{\ast}}

\newcommand{\Zs}{Z^{\ast}}
\newcommand{\zs}{z^{\ast}}
\newcommand{\eps}{\varepsilon}

\calclayout

\begin{document}

\title[]{Diametral strong diameter two property\\ of Banach spaces is stable under\\ direct sums with 1-norm}%
\author{Rainis Haller}%
\address{Institute of Mathematics and Statistics, University of Tartu,\newline J.~Liivi 2, 50409 Tartu, Estonia}%
\email{\tiny rainis.haller@ut.ee, katriinp@ut.ee, mart.poldvere@ut.ee}%
\thanks{The research was supported by institutional research funding IUT20-57
of the Estonian Ministry of Education and Research.}

\author{Katriin Pirk}%
%\address{}%
%\email{katriinp@ut.ee}%

\author{M\"{a}rt P\~{o}ldvere}%
%\email{mart.poldvere@ut.ee}%

%\thanks{}%

%\date{}%
%\dedicatory{}%
%\commby{}%
% ----------------------------------------------------------------
\begin{abstract}
We prove that the diametral strong diameter $2$ property of a Banach space (meaning that, in convex combinations of relatively weakly open subsets of its unit ball, every point has an ``almost diametral'' point) is stable under $1$-sums, i.e., the direct sum of two spaces with the diametral strong diameter $2$ property equipped with the $1$-norm has again this property.
\end{abstract}
\maketitle
% ----------------------------------------------------------------
%\section{}
All Banach spaces considered in this note are over the real field.
The closed unit ball and the unit sphere of a Banach space~$X$ will be denoted by $B_X$ and $S_X$, respectively.
Whenever referring to a relative weak topology, we mean such topology on the closed unit ball of the space under consideration.

Different versions of \emph{diameter 2 properties} for a Banach space mean that certain subsets of its unit ball (e.g., slices, nonempty relatively weakly open subsets, or convex combinations of weakly open subsets) have diameter equal to $2$. In recent years, these properties have been intensively studied (see, e.g., \cite{NW, P, ALN, 1509.02061, MR3334951, HL, MR3346197, ABL, MR3281132, BGLPRZ, AHNTT} for some typical results and further references).

To clarify the cap between the well-studied Daugavet property \cite{DW2} and known diameter $2$ properties,
the \emph{diametral} diameter $2$ properties
%which are strengthened versions of diameter $2$ properties in terms of abundance of diametral points
were introduced and studied in the recent preprint \cite{1509.02061}.
In particular, the stability under $p$-sums of diametral diameter $2$ properties was analyzed.
The question whether the $1$-sum of two Banach spaces enjoying the diametral \emph{strong} diameter $2$ property also has this property, was posed as an open problem  in \cite{1509.02061}.
Below, we shall answer this question in the affirmative.
%but the case $p=1$ remained open for the diametral \emph{strong} diameter $2$ property.

\begin{definition}[\cite{1509.02061}]
A Banach space $X$ is said to have the \emph{diametral diameter $2$ property} (briefly, DSD2P)
if, given $n\in\N$, relatively weakly open subsets $U_1,\dotsc,U_n$ of $B_X$,
$\lambda_1,\dotsc,\lambda_n\in[0,1]$ with $\sum_{i=1}^n{\lambda_i}=1$, $x\in\sum_{i=1}^n\lambda_i U_i$, and $\eps>0$,
there is a $u\in\sum_{i=1}^n \lambda_i U_i$ satisfying
\[
\|x-u\|\geq\|x\|+1-\eps.
\]
\end{definition}

\begin{theorem}
Suppose that Banach spaces $X$ and $Y$ have the DSD2P.
Then also the $1$-sum $X\oplus_1 Y$ has the DSD2P.
\end{theorem}

Our proof of Theorem makes use of the following observation:
\begin{itemize}
\begin{comment}
\item[(1)]
the space $X$ is infinite dimensional (because clearly no finite dimensional space can have the DSD$2$P);
%
\item[(2)]
the sets $U_1,\dotsc,U_n$ are convex;
%
\end{comment}
\item[($\bullet$)]
in Definition, one may assume that the element $x$ is of the form  $x=\sum_{i=1}^n\lambda_i x_i$ where $x_i\in S_X\cap U_i$.
\end{itemize}

\noindent
For ($\bullet$), first notice that the space $X$ may be assumed to be infinite dimensional (because clearly no finite dimensional space can have the DSD$2$P) and the sets $U_1,\dotsc,U_n$ to be convex. Now,
for ($\bullet$), it suffices to observe that
\begin{itemize}
\item[$(\circ)$]
every $a\in U_i$ can be written in the form $a=(1-\mu_i) y_i + \mu_i z_i$ where $\mu_i\in[0,1]$ and $y_i,z_i\in S_X\cap U_i$,
\end{itemize}
\noindent because, if $(\circ$) holds, then the element $x$ can be written as
\[
x=\sum_{i=1}^n\lambda_i(1-\mu_i)y_{i}+\sum_{i=1}^n\lambda_i\mu_i z_{i}
\]
and (by the convexity of $U_1,\dotsc,U_n$)
\[
\sum_{i=1}^n\lambda_i(1-\mu_i)U_i+\sum_{i=1}^n\lambda_i\mu_iU_i\subset\sum_{i=1}^n\lambda_i U_i.
\]

It remains to prove $(\circ)$. Let $i\in\{1,\dotsc,n\}$ and let $a\in U_i$, $\|a\|<1$.
Let $m\in\N$, $\xs_{1},\dotsc,\xs_{m}\in\Xs$, and $\delta>0$ be such that
\[
U_i\supset\bigl\{b\in B_X\colon\,|\xs_{j}(b)-\xs_{j}(a)|<\delta,\,j=1,\dotsc,m\bigr\}.
\]
Choose a non-zero $c\in \bigcap_{j=1}^{m}\ker\xs_{j}$ (such $c$ exists when the space $X$ is infinite dimensional),
and consider the function $f(t)=\|a+tc\|$, $t\in\R$. Since $f(0)=\|a\|<1$ and $f(t)\xrightarrow[t\to\pm\infty]{}\infty$, there are $s,t\in(0,\infty)$ such that $f(-s)=f(t)=1$,
but now $y_{i}:=a-sc$, $z_{i}:=a+tc$, and $\mu_i:=\dfrac{s}{s+t}$ do the job.

\begin{proof}[Proof of Theorem]
Put $Z:=X\oplus_1 Y$, and let $n\in\N$, let $W_1,\dotsc W_n$ be relatively weakly open subsets of $B_Z$,
let $\lambda_1,\dotsc,\lambda_n\in[0,1]$ satisfy $\sum_{i=1}^n{\lambda_i}=1$,
and let $z=\sum_{i=1}^n\lambda_i z_i$ where $z_i=(x_i,y_i)\in S_Z\cap W_i$.
We must find a $w=(u,v)\in\sum_{i=1}^n\lambda_i W_i$ so that $\|z-w\|\geq\|z\|+1-\eps$,
i.e., putting $x:=\sum_{i=1}^n\lambda_i x_i$ and $y:=\sum_{i=1}^n\lambda_i y_i$ (now one has $z=(x,y)$),
\[
\|x-u\|+\|y-v\|\geq\|x\|+\|y\|+1-\eps.
\]

For every $i\in\{1,\dotsc,n\}$, putting
\[
\widehat{x}_i=
\begin{cases}
\dfrac{x_i}{\|x_i\|},&\quad\text{if $x_i\not=0$,}\\
0,&\quad\text{if $x_i=0$,}
\end{cases}
\qquad\text{and}\qquad
\widehat{y}_i=
\begin{cases}
\dfrac{y_i}{\|y_i\|},&\quad\text{if $y_i\not=0$,}\\
0,&\quad\text{if $y_i=0$,}
\end{cases}
\]
there are relatively weakly open neighbourhoods $U_i\subset B_X$ and $V_i\subset B_Y$ of $\widehat{x}_i$ and $\widehat{y}_i$, respectively,
such that $\bigl(\|x_i\|\,U_i\bigr)\times\bigl(\|y_i\|\,V_i\bigr)\subset W_i$. %(we emphasize that $\|x_i\|\,U_i$ and $\|y_i\|\,V_i$ may fail to be relatively weakly open)
Indeed,
letting $m\in\N$, $\zs_j=(\xs_j,\ys_j)\in S_{\Zs}$, $j=1,\dotsc,m$, and $\delta>0$ be such that
\[
W_i\supset\bigl\{w\in B_Z\colon\,|\zs_j(w)-\zs_j(z_i)|<\delta,\,j=1,\dotsc,m\bigr\}
\]
and defining
\begin{align*}
U_i&:=\bigl\{u\in B_X\colon\,|\xs_j(u)-\xs_j(\widehat{x}_i)|<\delta,\,j=1,\dotsc,m\bigr\},\\
V_i&:=\bigl\{v\in B_Y\colon\,|\ys_j(v)-\ys_j(\widehat{y}_i)|<\delta,\,j=1,\dotsc,m\bigr\},
\end{align*}
one has, whenever $u\in U_i$ and $v\in V_i$, for every $j\in\{1,\dotsc,m\}$,
\begin{align*}
\Bigl|\zs_j\bigl(\|x_i\|u,&\|y_i\|v\bigr)-\zs_j(z_i)\Bigr|
=\Bigl|\zs_j\bigl(\|x_i\|u,\|y_i\|v\bigr)-\zs_j(x_i,y_i)\Bigr|\\
&=\Bigl|\xs_j\bigl(\|x_i\|u\bigr)+\ys_j\bigl(\|y_i\|v\bigr)-\xs_j(x_i)-\ys_j(y_i)\Bigr|\\
&=\Bigl|\xs_j\bigl(\|x_i\|u\bigr)+\ys_j\bigl(\|y_i\|v\bigr)-\xs_j\bigl(\|x_i\|\widehat{x}_i\bigr)-\ys_j\bigl(\|y_i\|\widehat{y}_i\bigr)\Bigr|\\
&=\bigl|\|x_i\|\xs_j(u-\widehat{x}_i)+\|y_i\|\ys_j(v-\widehat{y}_i)\bigr|\\
&\leq\|x_i\|\,\bigl|\xs_j(u-\widehat{x}_i)\bigr|+\|y_i\|\,\bigl|\ys_j(v-\widehat{y}_i)\bigr|\\
&<\bigl(\|x_i\|+\|y_i\|\bigr)\delta=\|z_i\|\delta\\
&=\delta.
\end{align*}

Put
\[
\alpha:=\sum_{i=1}^n\lambda_i\|x_i\|
\quad\text{and}\quad
\beta:=\sum_{i=1}^n\lambda_i\|y_i\|.
\]
Notice that
\[
\alpha+\beta=\sum_{i=1}^n\lambda_i\bigl(\|x_i\|+\|y_i\|\bigr)=\sum_{i=1}^n\lambda_i\|z_i\|=\sum_{i=1}^n\lambda_i=1.
\]
We only consider the case when both $\alpha\not=0$ and $\beta\not=0$. (The case when $\alpha=0$ or $\beta=0$ can be handled similarly and is, in fact, simpler.)

For every $i\in\{1,\dotsc,n\}$, letting
\[
\alpha_i:=\frac{\lambda_i\|x_i\|}{\alpha}
\quad\text{and}\quad
\beta_i:=\frac{\lambda_i\|y_i\|}{\beta},
\]
one has $\alpha_i,\beta_i\in[0,1]$, and $\sum_{i=1}^n\alpha_i=\sum_{i=1}^n\beta_i=1$.
Since $X$ and $Y$ have the DSD$2$P, observing that
\begin{align*}
\frac{x}{\alpha}=\sum_{i=1}^n\frac{\lambda_i\|x_i\|}{\alpha}\,\widehat{x}_i\in\sum_{i=1}^n\alpha_i U_i
\quad\text{and}\quad
\frac{y}{\beta}=\sum_{i=1}^n\frac{\lambda_i\|y_i\|}{\beta}\,\widehat{y}_i\in\sum_{i=1}^n\beta_i V_i,
\end{align*}
there are $u_0\in\sum_{i=1}^n\alpha_i U_i$ and $v_0\in\sum_{i=1}^n\beta_i V_i$ such that
\begin{align*}
\Bigl\|\dfrac{x}{\alpha}-u_0\Bigr\|\geq\frac{1}{\alpha}\|x\|+1-\eps
\quad\text{and}\quad
\Bigl\|\dfrac{y}{\beta}-v_0\Bigr\|\geq\frac{1}{\beta}\|y\|+1-\eps.
\end{align*}
Finally, putting
\begin{align*}
u&:=\alpha u_0\in\sum_{i=1}^n\alpha \alpha_i U_i=\sum_{i=1}^n\lambda_i\|x_i\| U_i,\\
v&:=\beta v_0\in\sum_{i=1}^n\beta \beta_i V_i=\sum_{i=1}^n\lambda_i\|y_i\| V_i,
\end{align*}
one has
\[
(u,v)\in\sum_{i=1}^n\lambda_i\Bigl(\bigl(\|x_i\|\,U_i\bigr)\times\bigl(\|y_i\|\,V_i\bigr)\Bigr)\subset\sum_{i=1}^n\lambda_i W_i
\]
and
\[
\|x-u\|+\|y-v\|\geq\|x\|+\|y\|+(\alpha+\beta)(1-\eps)=\|x\|+\|y\|+1-\eps,
\]
as desired.
\end{proof}

Thus the stability of the diametral strong diameter $2$ property under $1$- and $\infty$-sums is similar to that of the Daugavet property.
In fact, among all 1-unconditional sums of two Daugavet spaces only the $1$- and $\infty$-sum have the Daugavet property. Whether the diametral strong diameter two property and the Daugavet property coincide remains an open question.

\bibliographystyle{amsplain}
\footnotesize
%\bibliography{Bibliography}

{\providecommand{\noopsort}[1]{}}
\providecommand{\bysame}{\leavevmode\hbox to3em{\hrulefill}\thinspace}
\providecommand{\MR}{\relax\ifhmode\unskip\space\fi MR }
\providecommand{\MRhref}[2]{%
  \href{http://www.ams.org/mathscinet-getitem?mr=#1}{#2}
}
\providecommand{\href}[2]{#2}

\end{document}